\documentclass[MSNbibl,number,citesort,seceqn,dvips]{arxbj}
\usepackage{upgreek}
\usepackage{mathrsfs}

\aid{0}
\volume{17}
\issue{3}
\pubyear{2011}
\firstpage{916}
\lastpage{941}
\doi{10.3150/10-BEJ311}

\makeatletter

\newcommand{\R}{\mathbb{R}}
\newcommand{\morf}[4][\to]{ #2 \dvtx  #3 #1 #4}
\renewcommand{\d}{\mathrm{d}}
\newcommand{\var}[1]{\mathrm{V}_{#1}}
\newcommand{\ac}[1]{\mathrm{R}_{#1}}
\newcommand{\acc}[1]{\mathrm{\bar R}_{#1}}
\newcommand{\e}{\mathrm{e}}
\newcommand{\norm}[2][]{#1\Vert #2#1\Vert}
\newcommand{\abs}[2][]{#1\vert #2#1\vert}
\newcommand{\F}{\mathcal F}
\newcommand{\spn}{\overline{\operatorname{span}}}
\newcommand{\eqref}[1]{(\ref{#1})}
\newcommand{\vart}{\mathrm{V}_N''}
\newcommand{\varo}{\mathrm{V}_N'}
\newcommand{\varte}{\mathrm{V}_N'''}

\newtheorem{lemma}{Lemma}[section]
\newtheorem{proposition}[lemma]{Proposition}
\newtheorem{theorem}[lemma]{Theorem}
\newtheorem{corollary}[lemma]{Corollary}
\newremark{example}[lemma]{Example}

\newremark{remark}[lemma]{Remark}

\def\bptnote#1{}

\makeatother

\begin{document}
\begin{frontmatter}

\title{Quasi Ornstein--Uhlenbeck processes}

\runtitle{Quasi Ornstein--Uhlenbeck processes}

\begin{aug}
\author{\inits{O. E.}\fnms{Ole E.} \snm{Barndorff-Nielsen}\thanksref{e1}%
\ead[label=e1,mark]{oebn@imf.au.dk}}
\and
\author{\inits{A.}\fnms{Andreas} \snm{Basse-O'Connor}\corref{}\thanksref{e2}%
\ead[label=e2,mark]{basse@imf.au.dk}}

\runauthor{O.E. Barndorff-Nielsen and A. Basse-O'Connor}

\address{Department of Mathematical Sciences, University of Aarhus
and Thiele Centre, Ny Munkegade 118,
DK-8000 Aarhus C, Denmark. \printead{e1,e2}}

\end{aug}

\received{\smonth{12} \syear{2009}}
\revised{\smonth{6} \syear{2010}}

%
\begin{abstract}
The question of existence and properties of stationary solutions to
Langevin equations driven by noise processes with stationary increments
is discussed, with particular focus on noise processes of
pseudo-moving-average type. On account of the Wold--Karhunen
decomposition theorem, such solutions are, in principle, representable
as a moving average (plus a drift-like term) but the kernel in the
moving average is generally not available in explicit form. A class of
cases is determined where an explicit expression of the kernel can be
given, and this is used to obtain information on the asymptotic
behavior of the associated autocorrelation functions, both for small
and large lags. Applications to Gaussian- and L\'{e}vy-driven
fractional Ornstein--Uhlenbeck processes are presented. A Fubini
theorem for L\'{e}vy bases is established as an element in the
derivations.
\end{abstract}

%

\end{frontmatter}

\section{Introduction}

This paper studies the existence and properties of stationary solutions
to Langevin equations driven by a noise process $N$ with, in general,
stationary dependent increments. We shall refer to such solutions as
quasi Ornstein--Uhlenbeck (QOU) processes. Of particular interest are
the cases where the noise process is of the pseudo-moving-average (PMA)
type. In wide generality, the stationary solutions can, in principle, be
written in the form of a Wold--Karhunen-type representation, but it is
relatively rare that an explicit expression for the kernel of such a
representation can be given. When this is possible it often provides a
more direct and simpler access to the character and properties of the
process; for instance, concerning the autocovariance function.

This will be demonstrated in applications to the case where the noise
process $N$ is of the pseudo-moving-average kind, including fractional
Brownian motion and, more generally, fractional L\'{e}vy motions. Of
some particular interest for turbulence theory is the large and small
lags limit behavior of the autocovariance function of the
Ornstein--Uhlenbeck-type process driven by fractional Brownian motion,
which has been proposed as a representation of homogeneous Eulerian
turbulent velocities; see Shao \cite{Shao95}.

The fractional Brownian and L\'{e}vy motions are not of the
semimartingale type. Another
non-semimartingale case covered is $N_{t}=\int_{\mathcal X
}B_{t}^{(x)}m(\mathrm{d}x),$
where the processes $B_{\cdot}^{( x) }$ are Brownian
motions in different filtrations and $m$ is a measure on some space
$\mathcal X$.

In recent applications of stochastics, particularly in finance and in
turbulence, modifications of classic noise processes by time change or
by volatility modifications are of central importance; see
Barndorff-Nielsen and Shephard \cite{BNSBook} and Barndorff-Nielsen and Shiryaev \cite{BNShi10} and references given therein.
Prominent examples of such processes are $\mathrm{d}N_{t}=\sigma_{t}\,\mathrm{d}B_{t}$, where $B$ is Brownian motion and $\sigma>0$ is a
predictable stationary process -- for instance, the square root of a
superposition of inverse Gaussian Ornstein--Uhlenbeck processes (cf.
Barndorff-Nielsen and Shephard \cite{OEBNnonGOU} and Barndorff-Nielsen and Stelzer \cite{OEBNRob}) -- and $N_{t}=L_{T_{t}}$,
where $L$ is a L\'{e}vy process and $T$ is a time change process with
stationary increments (cf. Carr
  \textit{et~al}. \cite{CGMY03}). The theory discussed in
the present paper applies also to processes of this type.

The structure of the paper is as follows. Section~\ref{exi} defines the
concept of QOU processes and provides conditions for existence and
uniqueness of stationary solutions to the Langevin equation. The form
of the autocovariance function of the solutions is given and its
asymptotic behavior for \thinspace$t\rightarrow\infty$ is discussed.
As an intermediate step, a Fubini theorem for L\'{e}vy bases is
established in Section~\ref{fubini}. In Section~\ref{ma} explicit forms
of Wold--Karhunen representations are derived and used to analyze the
asymptotics, under more specialized assumptions, of the autocovariance
functions, both for $t\rightarrow\infty$ and for $t\rightarrow0$.
The \hyperref[app]{Appendix} establishes an auxiliary continuity result of a technical
nature.

\section{Langevin equations and QOU processes}\label{exi}

Let $N=(N_t)_{t\in\R}$ be a measurable process with stationary
increments and let $\lambda>0$ be a~positive number. By a QOU process
$X$ driven by $N$ and with parameter $\lambda$, we mean a~stationary
solution to the Langevin equation $\d{X} _t=-\lambda
X_t\,\d{t}+\d{N}_t$, that is, $X=(X_t)_{t\in\R}$ is a stationary
process that satisfies
%
\begin{equation}
\label{eq:5}
X_t=X_0-\lambda\int_0^t X_s\,\d{s} +N_t,\qquad t\in\R,
\end{equation}
where the integral is a pathwise Lebesgue integral.
For all $a<b$ we use the notation $\int_{b}^a :=-\int_a^b$.
Recall that a process $Z=(Z_t)_{t\in\R}$ is measurable if
$(t,\omega)\mapsto Z_t(\omega)$ is
$(\mathcal B(\R)\otimes\mathcal F,\break\mathcal B(\R))$-measurable, and
that $Z$ has stationary increments if, for all $s\in\R$,
$(Z_t-Z_0)_{t\in\R}$ has the same finite distributions as
$(Z_{t+s}-Z_s)_{t\in\R}$. For $p>0$ we will say that a process $Z$ has
finite $p$ moments if $\mathrm{E}[|Z_t|^p]<\infty$ for all $t\in\R$. Moreover,
for $t\to0$ or $\infty$, we will write $f(t)\sim g(t)$, $f(t)=\mathrm{o}(g(t))$
or $f(t)=\mathrm{O}(g(t)),$ provided that $f(t)/g(t)\rightarrow1$,
$f(t)/g(t)\to0$ or $\limsup_{t}\abs{f(t)/g(t)}<\infty$, respectively.
For each process $Z$ with finite second moments, let $\var{Z}(t)
=\operatorname{Var}(Z_t)$ denote its variance function. When $Z$, in
addition, is stationary, let $\ac{Z}(t)=\operatorname{Cov}(Z_t, Z_0)$ denote
its autocovariance function and
$\acc{X}(t)=\ac{X}(0)-\ac{X}(t)=\frac{1}{2}\mathrm{E}[(X_t-X_0)^2]$ its
complementary autocovariance function.

Before discussing the general setting further, we recall some
well-known cases. The stationary solution $X$ to \eqref{eq:5} when
$N_t= \mu t+\sigma B_t$ (with $B$ the Brownian motion) is the Gaussian
Ornstein--Uhlenbeck process, $\mu/\lambda$ is the mean level,
$\lambda$
is the speed of reversion and $\sigma$ is the volatility. When $N$ is
a L\'{e}vy process, the corresponding QOU process, $X$, exists if and
only if $\mathrm{E}[\log^+|N_1|]<\infty$ or, equivalently, if and only if
$\int_{\{|x|>1\}}\log|x|\nu(\d{x})<\infty,$ where $\nu$ is the
L\'{e}vy measure of $N$; see Sato and Yamazato \cite{SatoOU} or Wolfe \cite{WOlfe}. In this
case $X$ is called an Ornstein--Uhlenbeck-type process; for
applications of such processes in financial economics, see
Barndorff-Nielsen and Shephard \cite{OEBNnonGOU,BNSBook}.

\subsection{Existence and uniqueness of QOU processes}

The first result below shows the existence and uniqueness for the
stationary solution $X$ to the Langevin equation $\d{X} _t=-\lambda
X_t\,\d{t}+\d{N}_t$ in the case where the noise $N$ is integrable. That
is, we show existence and uniqueness of QOU processes $X$. Moreover,
we provide an explicit form of the solution that is used to calculate
the mean and variance of $X$.

\begin{theorem}\label{lemma_Masani}
Let $N$ be a measurable process with stationary
increments and finite first moments, and let
$\lambda>0$ be a positive real number. Then, $X=(X_t)_{t\in\R}$,
given by
%
\begin{equation}
\label{eq:88}
X_t=N_t-\lambda \e^{-\lambda t}\int_{-\infty}^t \e^{\lambda s}
N_s\,\d{s},\qquad t\in\R,
\end{equation}
is
a QOU process driven by $N$ with
parameter $\lambda$ (the integral is a pathwise Lebesgue
integral). Furthermore, any other QOU process driven by
$N$ and with parameter $\lambda$ equals~$X$ in law. Finally,
if $N$ has finite $p$ moments, $p\geq1$, then
$X$ also has finite $p$ moments and is continuous in $L^p$.
\end{theorem}

\begin{remark}
It is an open problem to relax the integrability of $N$ in
Theorem~\ref{lemma_Masani}; that is, is it enough that $N$ has finite
$\log$ moments? Recall that when $N$ is a L\'{e}vy process, finite
$\log$ moments is a necessary and sufficient
condition for the existence of the corresponding
Ornstein--Uhlenbeck-type process.
\end{remark}

\begin{pf*}{Proof of Theorem \ref{lemma_Masani}}
Existence: Let $p\geq1$ and assume that $N$ has finite $p$ moments.
Choose $\alpha,\beta>0$, according to Corollary~\ref{p_cont}, such that
$\norm{N_t}_p\leq\alpha+\beta|t|$ for all $t\in\R$. By
Jensen's inequality,
\begin{eqnarray*}
\mathrm{E}\biggl[\biggl(\int_{-\infty}^t \e^{\lambda s}|N_s|\,\d{s}\biggr)^p\biggr]
&\leq&
(\e^{\lambda t}/\lambda)^{p-1}\int_{-\infty}^t \e^{\lambda s}\mathrm{E}[|N_s|^p]\,\d{s}\\
&\leq&
(\e^{\lambda t}/\lambda)^{p-1}\int_{-\infty}^t \e^{\lambda s}(\alpha+\beta|s|)^p\,\d{s}<\infty,
\end{eqnarray*}
which shows that the integral in \eqref{eq:88}
exists almost surely as a Lebesgue
integral and that $X_t$, given by \eqref{eq:88}, is
$p$-integrable. Using substitution we obtain from
\eqref{eq:88},
%
\begin{equation}
\label{eq:23}
X_t=\lambda\int_{-\infty}^0 \e^{\lambda
u}(N_t-N_{t+u})\,\d{u},\qquad t\in\R.
\end{equation}
By Corollary~\ref{p_cont}, $N$ is
$L^p$-continuous and, therefore, it follows that
the right-hand side of \eqref{eq:23} exists as a limit of
Riemann sums in $L^p$. Hence the
stationarity of the increments of $N$ implies that $X$ is
stationary. Moreover,
using integration by parts on $t\mapsto\int_{-\infty}^t \e^{\lambda s}
N_s(\omega)\,\d{s}$, we get
\[
\int_0^t X_s\,\d{s}=\e^{-\lambda t}\int_{-\infty}^t \e^{\lambda
s}N_s\,\d{s}-\int_{-\infty}^0 \e^{\lambda s}N_s\,\d{s},
\]
which shows that $X$ satisfies \eqref{eq:5},
and hence $X$ is a QOU process driven by $N$ with
parameter~$\lambda$.

Since $X$ is a measurable process with stationary increments and
finite $p$ moments, Proposition~\ref{p_cont} shows that it is
continuous in $L^p$.

To show uniqueness in law, let $Y$ be a QOU process
driven by $N$ with parameter $\lambda>0$, that is, $Y$ is a
stationary process that satisfies
\eqref{eq:5}. For all $t_0\in\R$
we have, with $Z_t=N_t-N_{t_0}+Y_{t_0}$, that
%
\begin{equation}
\label{eq:48}
Y_t=Z_t-\lambda\int_{t_0}^t Y_s\,\d{s},\qquad t\geq t_0.
\end{equation}
Solving \eqref{eq:48} pathwise, it follows that for all $t\geq t_0$,
\begin{eqnarray*}
Y_t&=& Z_t-\lambda \e^{-\lambda t}\int_{t_0}^t \e^{\lambda s} Z_s\,\d{s}
\\
&=& N_t-\lambda \e^{-\lambda t}\int_{t_0}^t \e^{\lambda s} N_s\,\d{s}+(Y_{t_0}-N_{t_0})\e^{-\lambda(t-t_0)} .
\end{eqnarray*}
Note that $\lim_{t\to\infty}(Y_{t_0}-N_{t_0})\e^{-\lambda
(t-t_0)}=0$ a.s., thus for all $n\geq1$ and $t_0<t_1<\cdots<t_n$, the
stationarity of $Y$ implies that for $k\to\infty$,
$(Y_{t_i+k})_{i=1}^n\Rightarrow(Y_{t_i})_{i=1}^n$ (for all random
vectors, $\Rightarrow$ denotes convergence in distribution). Therefore,
as $k\to\infty$,
\[
\biggl( N_{t_i+k}-\lambda \e^{-\lambda(t_i+k)}\int_{t_0}^{t_i+k}
\e^{\lambda s} N_s\,\d{s}\biggr)_{i=1}^n
\Rightarrow(Y_{t_i})_{i=1}^n.
\]
This shows that the distribution of $Y$ only depends
on $N$ and $\lambda$, and completes the proof.
\end{pf*}

Proposition~2.1 in Surgailis \textit{et~al}. \cite{Rosinskipreprint} and Proposition~2.1 in
Maejima and Yamamoto \cite{MaejimafsOU} provide existence results for stationary solutions
to Langevin equations.
However, these results do not cover Theorem~\ref{lemma_Masani}.
The first result considers only Bochner-type integrals and
the second result requires, in particular, that the sample paths of
$N$ are Riemann integrable.

Let $B=(B_t)_{t\in\R}$ denote an $\F$-Brownian motion indexed by $\R$
and $\sigma=(\sigma_t)_{t\in\R}$ be a~predictable process; that is,
$\sigma$ is measurable with respect to
\[
\mathscr P=\sigma\bigl((s,t]\times A\dvt s,t\in\R,s<t,A\in\F_s\bigr).
\]
Assume that, for all $u\in\R$, $(\sigma_t,B_t)_{t\in\R}$ has the same
finite-dimensional distributions as
$(\sigma_{t+u},B_{t+u}-B_u)_{t\in\R}$ and that $\sigma_0\in L^2$. Then
$N$, given by
%
\begin{equation}
\label{eq:125}
N_t=\int_0^t \sigma_s\,\d{B}_s,\qquad t\in\R,
\end{equation}
is a well-defined continuous process
with stationary increments and finite second moments. (Recall that for
$t<0$, $\int_0^t:=-\int_{t}^0$.)

\begin{corollary}\label{cor_sigma}
Let $N$ be given by \eqref{eq:125}. Then,
there exists a unique-in-law QOU process $X$ driven by $N$ with
parameter $\lambda>0$, and $X$ is given by
%
\begin{equation}
\label{eq:127}
X_t=\int_{-\infty}^t \e^{-\lambda(t-s)}\sigma_s\,\d{B}_s,\qquad
t\in\R.
\end{equation}
\end{corollary}

\begin{pf}
Since $N$ is a measurable process with finite second moments, it
follows by Theorem~\ref{lemma_Masani}
that there exists a unique-in-law QOU process $X$, and it is
given by
\begin{eqnarray}\label{eq:128}
X_t
&=&
N_t-\lambda \e^{-\lambda t}\int_{-\infty}^t \e^{\lambda s}N_s\,\d{s} =\lambda\int_{-\infty}^0 \e^{\lambda s}(N_t-N_{t+s})\,\d{s}\nonumber
\\[-8pt]\\[-8pt]
&=&
\lambda\int_{-\infty}^0 \biggl(\int_\R1_{(t+s,t]}(u) \e^{\lambda s}\sigma_u\,\d{B}_u\biggr)\,\d{s}.\nonumber
\end{eqnarray}
By an extension of
Protter \cite{Protter}, Chapter~IV, Theorem~65, from finite intervals to
infinite intervals we may switch the order of
integration in \eqref{eq:128} and hence we obtain \eqref{eq:127}.
\end{pf}

Let us conclude this section with formulas for the mean and variance of
a QOU process~$X$. In the rest of this section, let $N$ be a measurable
process with stationary increments and finite first moments and let $X$
be a QOU process driven by $N$ with parameter $\lambda>0$ (which exists
by Theorem~\ref{lemma_Masani}). Since $X$ is unique in law, it makes
sense to consider the mean and variance function of $X$. Let us assume
for simplicity that $N_0=0$ a.s. The following proposition gives the
mean and variance of $X$.

\begin{proposition}\label{m_and_v}
Let $N$ and $X$ be given as above. Then,
\[
\mathrm{E}[X_0]= \frac{\mathrm{E}[N_1]}{\lambda} \quad\mbox{and}\quad
\operatorname{Var}(X_0)=\frac{\lambda}{2}\int_0^\infty \e^{-\lambda
s}\var{N}(s)\,\d{s}.
\]
In the part concerning the variance of $X_0$, we assume that $N$ has
finite second moments.
\end{proposition}

Note that Proposition~\ref{m_and_v} shows that the variance of $X_0$ is
$\lambda/2$ times the Laplace transform of $\var{N}$.
In particular, if $N_t=\mu t+\sigma B^{H}_t$, where
$B^{H}$ is a fractional
Brownian motion (fBm) of index $H\in(0,1)$ (see
\cite{StoccalfBm} or \cite{Nua06}
for properties of the fBm), then $\mathrm{E}[N_1]=\mu$ and
$\var{N}(s)=\sigma^2\abs{s}^{2H}$ and hence, by
Proposition~\ref{m_and_v}, we have that
%
\begin{equation}
\label{eq:119}
\mathrm{E}[X_0]=\frac{\mu}{\lambda} \quad\mbox{and}\quad
\operatorname{Var}(X_0)=\frac{\sigma^2\Gamma(1+2H)}{2\lambda^{2H}}.
\end{equation}
For $H=1/2$, \eqref{eq:119} is well known, and in this case
$\operatorname{Var}(X_0)=\sigma^2/(2\lambda)$.

Before proving Proposition~\ref{m_and_v}, let us
note that $\mathrm{E}[N_t]=\mathrm{E}[N_1]t$ for all
$t\in\R$. Indeed, this follows by the continuity of $t\mapsto\mathrm{E}[N_t]$
(see Corollary~\ref{p_cont}) and the stationarity of the increments of
$N$.

\begin{pf*}{Proof of Proposition \ref{m_and_v}}
Recall that, by Corollary~\ref{p_cont}, we have that
$\mathrm{E}[\abs{N_t}]\leq\alpha+\beta\abs{t}$ for some $\alpha,\beta>0$.
Hence, by \eqref{eq:88} and Fubini's theorem, we have that
\begin{eqnarray*}
\mathrm{E}[X_0]
&=&
\mathrm{E}\biggl[-\lambda\int_{-\infty}^0 \e^{\lambda s} N_s\,\d{s}\biggr]=-\lambda\int_{-\infty}^0 \e^{\lambda s}\mathrm{E}[N_s]\,\d{s}
\\
&=&
-\lambda\mathrm{E}[N_1]\int_{-\infty}^0 \e^{\lambda s}s\,\d{s}=\mathrm{E}[N_1]/\lambda.
\end{eqnarray*}
This shows the part concerning the mean of
$X_0$.

To show the last part, assume that $N$ has finite second moments. By
using $\mathrm{E}[X_0]=\mathrm{E}[N_1]/\lambda$, \eqref{eq:88} shows that,
with $\tilde N_t:=N_t-\mathrm{E}[N_1]t$, we have
\[
\operatorname{Var}(X_0)=\mathrm{E}\bigl[(X_0-\mathrm{E}[X_0])^2\bigr]=
\mathrm{E}\biggl[\biggl(\lambda\int_{-\infty}^0 \e^{\lambda s} \tilde N_s\,\d{s}\biggr)^2 \biggr].
\]
Since $\norm{\tilde N_t}_2\leq
\alpha+\beta\abs{t}$ for some $\alpha,\beta>0$ by
Corollary~\ref{p_cont}, Fubini's theorem shows
\[
\operatorname{Var}(X_0)=\lambda^2\int_{-\infty}^0\int_{-\infty}^0
(\e^{\lambda
s}\e^{\lambda u} \mathrm{E}[\tilde N_s \tilde N_u])\,\d{s} \,\d{u},
\]
and since $\mathrm{E}[\tilde N_s \tilde
N_u]=\frac{1}{2}[\var{N}(s)+\var{N}(u)-\var{N}(s-u)],$ we have
\begin{eqnarray}\label{eq:121}
\operatorname{Var}(X_0)
&=&
\frac{\lambda^2}{2}\int_{-\infty}^0\int_{-\infty}^0
\bigl(\e^{\lambda
s}\e^{\lambda u} \bigl(\var{N}(s)+\var{N}(u)-\var{N}(s-u)\bigr)\bigr)\,\d{s}
\,\d{u}\nonumber
\\[-8pt]\\[-8pt]
&=&
\lambda\int_{-\infty}^0 \e^{\lambda
s}\var{N}(s)\,\d{s}-\frac{\lambda^2}{2}\int_{-\infty}^0
\e^{\lambda
u}\biggl(\int_{-\infty}^{-u}
\e^{\lambda(s+u)} \var{N}(s)\,\d{s}\biggr)\,\d{u}.\nonumber
\end{eqnarray}
Moreover,
\begin{eqnarray*}
&&
\frac{\lambda^2}{2}\int_{-\infty}^0 \e^{\lambda u}\biggl(\int_{-\infty}^{-u}\e^{\lambda(s+u)} \var{N}(s)\,\d{s}\biggr)\,\d{u}
\\
&&\quad=
\frac{\lambda^2}{2}\int_{\R} \var{N}(s)\e^{\lambda s}\biggl(\int_{-\infty}^{(-s)\wedge0}\e^{2\lambda u} \,\d{u}\biggr)\,\d{s}
\\
&&\quad=
\frac{\lambda^2}{2}\biggl(\int_{-\infty}^0 \var{N}(s)\e^{\lambda s}\biggl(\int_{-\infty}^{0} \e^{2\lambda u} \,\d{u}\biggr)\,\d{s}+\int_{0}^\infty\var{N}(s)\e^{\lambda
s}\biggl(\int_{-\infty}^{-s}\e^{2\lambda u} \,\d{u}\biggr)\,\d{s}\biggr)
\\
&&\quad=
\frac{\lambda}{4}\biggl(\int_{-\infty}^0 \var{N}(s)\e^{\lambda s}\,\d{s}+\int_{0}^\infty\var{N}(s)\e^{\lambda s}(\e^{-2\lambda s})\,\d{s}\biggr)
\\
&&\quad=
\frac{\lambda}{2}\int_0^\infty \e^{-\lambda s}\var{N}(s)\,\d{s},
\end{eqnarray*}
which, by \eqref{eq:121}, gives the expression for
the variance of $X_0$.
\end{pf*}

\subsection{Asymptotic behavior of the autocovariance function}\label{asy_auto}

The next result shows that the autocovariance function of a
QOU process $X$ driven by $N$ with parameter $\lambda$ has the same
asymptotic behavior at infinity as the second derivative of the
variance function of $N$ divided by $2\lambda^2$.

\begin{proposition}\label{var}
Let $N$ be a measurable process with stationary
increments, $N_0=0$ a.s., and finite second moments. Let $X$ be a QOU
process driven by $N$ with parameter $\lambda>0$.
\begin{longlist}[(ii)]
\item\hypertarget{infty}
Assume that $\var{N}$ is three times continuous differentiable
in a neighborhood of $\infty$,
and for $t\to \infty$ we have that $\vart(t)=\mathrm{O}(\e^{(\lambda/2)t})$,
$\e^{-\lambda t}=\mathrm{o}(\vart(t))$ and $\varte(t)=\mathrm{o}(\vart(t))$.
Then, for $t\rightarrow\infty$, we have $\ac{X}(t)\sim
(\frac{1}{2\lambda^2})\vart(t)$.
\item\hypertarget{zero} Assume for $t\to0$ that $t^2=\mathrm{o}(\var{N}(t))$,
then, for $t\to0,$ we have $\acc{X}(t)\sim\frac{1}{2} \var{N}(t)$.
More generally, let $p\geq1$ and assume that $N$ has finite $p$
moments and
$t=\mathrm{o}(\|N_t\|_p)$ as $t\to0$. Then, for $t\rightarrow0$, we have
$\|X_t-X_0\|_p\sim\|N_t\|_p$.
\end{longlist}
\end{proposition}

Note that by Proposition~\ref{var}\hyperlink{zero}{(ii)} the short-term
asymptotic behavior of $\acc{X}$ is not influenced by~$\lambda$.

\begin{pf*}{Proof of Proposition \ref{var}}
\hyperlink{infty}{(i)} Let $\beta>0$ and assume that $\var{N}$ is three times
continuous differentiable on $(\beta,\infty)$; that is, $\var{N}\in
C^3((\beta,\infty);\R)$. Let $t_0=\beta+1$, and let
us show that for $t\geq t_0$ and $t\to
\infty$,
%
\begin{equation}
\label{eq:100}
\ac{X}(t)=\frac{\e^{-\lambda t}}{4\lambda} \int_{t_0}^t \e^{\lambda u}
\vart(u)\,\d{u}
+\frac{\e^{\lambda t}}{4\lambda} \int_t^\infty \e^{-\lambda u} \vart
(u)\,\d{u}+\mathrm{O}(\e^{-\lambda t}).
\end{equation}
If we have shown \eqref{eq:100}, then, by using that
$\e^{-\lambda t}=\mathrm{o}(\vart(t))$, $\varte(t)=\mathrm{o}(\vart(t))$ and
l'H\^{o}pital's rule, \hyperlink{infty}{(i)} follows.

Similar to the proof of Proposition~\ref{m_and_v}, let $\tilde
N_t=N_t-\mathrm{E}[N_1]t$. To show \eqref{eq:100}, recall that by
Corollary~\ref{p_cont} we have $\norm{\tilde N_t}_2\leq\alpha+\beta
\abs{t}$ for some $\alpha,\beta>0$. Hence, by \eqref{eq:88} and
Fubini's theorem, we find that
%
\begin{equation}
\label{eq:66}
\ac{X}(t)=\mathrm{E}\bigl[(X_t-\mathrm{E}[X_t])(X_0-\mathrm{E}[X_0])\bigr]
=g(t)-\lambda \e^{-\lambda t}\int_{-\infty}^t
\e^{\lambda s}g(s)\,\d{s},
\end{equation}
where
\[
g(t)=-\lambda\int_{-\infty}^0 \e^{\lambda s}
\mathrm{E}[\tilde N_s \tilde N_t]\,\d{s},\qquad t\in\R.
\]
Since $\mathrm{E}[\tilde N_s
\tilde N_t]=\frac{1}{2}[\var{N}(t)+\var{N}(s)-\var{N}(s-t)],$ we have
that
\begin{eqnarray}\label{eq:68}
g(t)
&=&
-\frac{\lambda}{2} \int_{-\infty}^0 \e^{\lambda s} [\var{N}(t)+\var{N}(s)-\var{N}(t-s)]\,\d{s}\nonumber
\\[-8pt]\\[-8pt]
&=&
-\frac{1}{2}\biggl(\var{N}(t)-\lambda \e^{\lambda t}\int_{t}^\infty \e^{-\lambda s}\var{N}(s)\,\d{s}\biggr)\nonumber
-\frac{\lambda}{2}\int_{-\infty}^0 \e^{\lambda s}\var{N}(s)\,\d{s}.
\end{eqnarray}
From \eqref{eq:68}
it follows that $g\in C^1((\beta,\infty);\R)$ and hence, using partial
integration on \eqref{eq:66}, we have for $t\geq t_0$,
%
\begin{equation}
\label{eq:109}
\ac{X}(t)=\e^{-\lambda t}\int_{t_0}^t \e^{\lambda s}g'(s)\,\d
{s}+\e^{-\lambda t}\biggl(\e^{\lambda t_0}g(t_0)-\lambda
\int_{-\infty}^{t_0} \e^{\lambda s} g(s)\,\d{s}\biggr).
\end{equation}
Moreover, from \eqref{eq:68} and for $t\geq t_0$, we find
%
\begin{equation}
\label{eq:57}
g'(t)=-\frac{1}{2}\biggl(\mathrm{V'} _N(t)-\lambda^2\e^{\lambda
t}\int_t^\infty \e^{-\lambda s}\var{N}(s)\,\d{s}+\lambda
\var{N}(t)\biggr).
\end{equation}
For $t\to\infty$ we have, by assumption, that
$\vart(t)=\mathrm{O}(\e^{(\lambda/2)t})$, and hence $\varo
(t)=\mathrm{O}(\e^{(\lambda/2)t})$. Thus, from \eqref{eq:57} and a double use of
partial integration, we obtain that
%
\begin{equation}
\label{eq:70}
g'(t)
=\frac{\e^{\lambda t} }{2}\int_t^\infty \e^{-\lambda
s} \mathrm{V''}_N(s)\,\d{s},\qquad t\geq t_0.
\end{equation}
Using \eqref{eq:70}, Fubini's theorem and that $\vart
(t)=\mathrm{O}(\e^{(\lambda/2)t}),$ we have for $t\geq t_0$,
\begin{eqnarray*}
\label{eq:74}
&&
\e^{-\lambda t}\int_{t_0}^t \e^{\lambda s}g'(s)\,\d{s}
\\
&&\quad=
\e^{-\lambda t}\int_{t_0}^t \e^{\lambda s}\biggl(\frac{\e^{\lambda s} }{2}\int_s^\infty \e^{-\lambda u}\vart(u)\,\d{u}\biggr)\,\d{s}
\\
&&\quad=
\e^{-\lambda t} \int_{t_0}^\infty \e^{-\lambda u}
\vart(u)\biggl(\int_{t_0}^{t\wedge u} \frac{1}{2}\e^{2\lambda s}\,\d{s}\biggr)\,\d
{u}
\\
&&\quad=
\e^{-\lambda t} \int_{t_0}^\infty \e^{-\lambda u}
\vart(u)\biggl( \frac{1}{4\lambda}\bigl(\e^{2\lambda(t\wedge
u)}-\e^{2\lambda t_0}\bigr)
\biggr)\,\d{u}
\\
&&\quad=
\frac{\e^{-\lambda t}}{4\lambda} \int_{t_0}^t \e^{\lambda u}
\vart(u)\,\d{u}
+\frac{\e^{\lambda t}}{4\lambda} \int_t^\infty \e^{-\lambda u} \vart
(u)\,\d{u}-\e^{-\lambda t}\biggl( \frac{\e^{2\lambda
t_0}}{4\lambda}\int_{t_0}^\infty \e^{-\lambda u} \vart
(u)\,\d{u}\biggr).
\end{eqnarray*}
Combining this with \eqref{eq:109} we obtain
\eqref{eq:100}, and the proof of \hyperlink{infty}{(i)} is complete.

\hyperlink{zero}{(ii)} Using \eqref{eq:5} we have for all $t>0$ that
\[
\|X_t-X_0\|_p\leq\|N_t\|_p+\lambda\int_0^t \|X_s\|_p\,\d{s}=\|N_t\|_p+\lambda t\|X_0\|_p.
\]
On the other hand,
\[
\|X_t-X_0\|_p\geq \|N_t\|_p-\lambda\int_0^t\|X_s\|_p\,\d{s}=\|N_t\|_p-\lambda t\|X_0\|_p,
\]
which shows that
\[
1-\lambda\|X_0\|_p \frac{t}{\|N_t\|_p}\leq\frac{\|X_t-X_0\|_p}{\|N_t\|_p}\leq1+\lambda\|X_0\|_p \frac{t}{\|N_t\|_p}.
\]
A similar inequality is available when $t<0$,
and hence for $t\to0$ we have that $\|X_t-X_0\|_p\sim\|N_t\|_p$ if
$\lim_{t\to 0}(t/\|N_t\|_p)=0$.
\end{pf*}

When $N$ is an fBm of index $H\in(0,1),$ then $\var{N}(t)=|t|^{2H}$,
and hence
\[
\vart(t)=2H(2H-1)t^{2H-2},\qquad t>0.
\]
The conditions in Proposition~\ref{var} are
clearly fulfilled and thus we have the following corollary.

\begin{corollary}\label{cor_var}
Let $N$ be an fBm of index $H\in(0,1)$,
and let $X$ be a QOU process driven by $N$ with parameter
$\lambda>0$. For $H\in
(0,1)\setminus\{\frac{1}{2}\}$ and $t\rightarrow\infty$, we have
$\ac{X}(t)\sim(H(2H-1)/\lambda^2)t^{2H-2}$.
For $H\in(0,1)$ and $t\rightarrow0$,
we have $\acc{X}(t)\sim\frac{1}{2}|t|^{2H}$.
\end{corollary}

The above result concerning the behavior of $\ac{X}$ for $t\rightarrow
\infty$ when $N$ is an fBm has been obtained previously via a
different approach by Cheridito \textit{et~al.} \cite{CheKawMae03}, Theorem~2.3.

A square-integrable stationary process $Y=(Y_t)_{t\in\R}$ is said to
have \textit{long-range dependence} of order $\alpha\in(0,1)$
if $\ac{Y}$ is regularly varying at $\infty$ of index
$-\alpha$. Recall that a~function $\morf{f}{\R}{\R}$ is
regularly varying at $\infty$ of index
$\beta\in\R$ if, for $t\to\infty$, $f(t)\sim t^\beta l(t),$ where
$l$ is slowly varying, which means that for all $a>0$,
$\lim_{t\to\infty}l(at)/l(t)=1$. Many empirical observations have
shown evidence for long-range dependence in various fields, such as
finance, telecommunication and hydrology; see Doukhan \textit{et~al}. \cite{longrange}.
Let $X$ be a QOU process driven by $N$; then
Proposition~\ref{var}\hyperlink{infty}{(i)} shows that $X$ has long-range
dependence of order $\alpha\in(0,1)$ if and only if $\mathrm{V}_N''$
is regularly varying at $\infty$ of order $-\alpha$. Furthermore,
Proposition~\ref{cor_MA_noise}\hyperlink{pma_1}{(i)} below shows how to
construct
QOUs with long-range dependence. More precisely, if $X$ is a QOU driven
by $N,$ where $N$ is given by \eqref{pma_noise}, and for some $\alpha
\in(0,1)$ and $t\to\infty$, $f'(t)\sim ct^{(\alpha-1)/2}$, then $X$
has long-range dependence of order $\alpha$. The example $f(t)=(\delta
\vee t)^{H-1/2}$, with $\delta\geq0$ and $H\in(\frac{1}{2},1)$ is
considered in Corollary~\ref{exp_fBm}
and it follows that the QOU process $X$ has long-range dependence
of order $2-2H$. Here $X$ is a fractional Ornstein--Uhlenbeck process
if $\delta=0$, and a semimartingale if and only if $\delta>0$. A
quite different type of semimartingale with long-range dependence is
obtained for $N=\sigma\bullet B$ with $\sigma$ and $B$ independent
and $\sigma^{2}$ being a supOU process with long-range dependence, cf.
Barndorff-Nielsen \cite{OEBNRob1}, Barndorff-Nielsen and Stelzer \cite{OEBNRob} and
Barndorff-Nielsen and Shephard \cite{BNSBook}. Hence, by considering
more general processes than the fractional type, we can easily construct
stationary processes with long-range dependence within the
semimartingale framework.

\section{A Fubini theorem for L\'{e}vy bases}\label{fubini}

Let $\Lambda=\{\Lambda(A)\dvt A\in\mathcal S\}$ denote a centered L\'{e}vy
basis on a non-empty space $S$ equipped with a $\delta$-ring
$\mathcal S$, see Rajput and Rosi{\'{n}}ski \cite{Rosinskispec}. (A L\'{e}vy basis is an
infinitely divisible, independently scattered random measure. Recall
also that a $\delta$-ring on $S$ is a family of subsets of $S$ that is closed
under union, countable intersection and set difference). As usual, we
assume that $\mathcal S$ is $\sigma$-finite, meaning that there exists
$(S_n)_{n\geq1}\subseteq\mathcal S$ such that $\bigcup_{n\geq1} S_n=S$.
All integrals $\int_S f(s)\Lambda(\d{s})$ will be defined in the
sense of Rajput and Rosi{\'{n}}ski \cite{Rosinskispec}. We can now find a measurable
parametrization of L\'{e}vy measures $\nu(\d{u},s)$ on $\R$, a
$\sigma$-finite measure $m$ on $S$ and a positive measurable
function $\sigma^2\dvtx S\rightarrow\R_+$,
such that for all $A\in\mathcal S$,
%
\begin{equation}
\label{eq:36}
\mathrm{E}\bigl[\e^{\mathrm{i}y\Lambda(A)}\bigr]=\exp\biggl(\int_A\biggl[-\sigma^2 (s)y^2/2+\int_\R(
\e^{\mathrm{i}yu}-1-\mathrm{i}yu)\nu(\d{u},s)\biggr] m(\d{s})\biggr),\qquad y\in\R,
\end{equation}
see \cite{Rosinskispec}. Let $\phi\dvtx \R\times S\mapsto\R$ be given by
\[
\phi(y,s)=y^2\sigma^2(s)+\int_\R\bigl[(uy)^2 1_{\{|uy|\leq
1\}}+(2|uy|-1)1_{\{|uy|>1\}}\bigr]\nu(\d{u},s),
\]
and for all measurable functions $\morf{g}{S}{\R}$
define
\[
\|g\|_\phi=\inf\biggl\{c>0\dvt\int_S \phi(c^{-1}g(s),s)m(\d{s})\leq
1\biggr\}\in[0,\infty].
\]
Moreover, let $L^\phi=L^\phi(S,\sigma(\mathcal S),m)$ denote
the \textit{Musielak--Orlicz space}
of measurable functions $g$ with
\[
\int_S \biggl[g(s)^2\sigma^2(s)+\int_\R\bigl( |ug(s)|^2\wedge
|ug(s)|\bigr)\nu(\d{u},s)\biggr]m(\d{s})<\infty,
\]
equipped with the Luxemburg norm
$\|g\|_\phi$. Note that $g\in L^\phi$ if and
only if $\|g\|_\phi<\infty$, since $\phi(2x,s)\leq C\phi(x,s)$ for some
$C>0$ and all $s\in S, x\in\R$.
We refer to Musielak \cite{Musielak} for the basic properties of
Musielak--Orlicz spaces. When $\sigma^2\equiv0$ and $g\in L^\phi$,
Theorem~2.1 in Marcus and Rosi{\'{n}}ski \cite{RosinskiL1norm}
shows that $\int_S g(s)\Lambda(\d{s})$ is well defined, integrable
and centered and
\[
c_1 \|g\|_\phi\leq \mathrm{E}\biggl[\biggl|\int_S g(s)\Lambda(\d{s})\biggr|\biggr]\leq c_2 \|g\|_\phi,
\]
and we may choose $c_1=1/8$ and $c_2=17/8$.
Hence for general $\sigma^2$ it is easily seen that for all $g\in
L^\phi$, $\int_S g(s)\Lambda(\d{s})$ is well defined, integrable and
centered and
%
\begin{equation}
\label{eq:25}
\mathrm{E}\biggl[\biggl|\int_S
g(s)\Lambda(\d{s})\biggr|\biggr]\leq2c_2 \|g\|_\phi.
\end{equation}

Let $T$ denote a complete separable metric space, and $Y=(Y_t)_{t\in
T}$ be given by
\[
Y_t=\int_S f(t,s)\Lambda(\d{s}),\qquad t\in T,
\]
for some measurable function $f(\cdot,\cdot)$ for which the integrals
are well defined. Then we can choose a measurable modification of $Y$.
Indeed, the existence of a measurable modification of $Y$ is equivalent
to measurability of $(t\in T)\mapsto(Y_t\in L^0)$ according to
Theorem~3 and the remark in Cohn \cite{Cohn}. Hence, since $f$ is
measurable, the maps $(t\in T)\mapsto
(\|f(t,\cdot)-g(\cdot)\|_\phi\in\R)$ for all $g\in L^\phi$ are
measurable. This shows that $(t\in\R)\mapsto(f(t,\cdot)\in L^\phi
)$ is
measurable since $L^\phi$ is a separable Banach space. Hence by
continuity of $(f(t,\cdot)\in L^\phi) \mapsto(Y_t\in L^0)$, see
Rajput and Rosi{\'{n}}ski \cite{Rosinskispec}, it follows that $(t\in T)\mapsto(Y_t\in L^0)$
is measurable.

Assume that $\mu$ is a $\sigma$-finite measure on a complete and
separable metric space $T$. Then we have the following stochastic
Fubini result extending Rosi{\'n}ski \cite{Rosinski}, Lemma~7.1; P\'{e}rez-Abreu and Rocha-Arteaga \cite{Arteaga},
Lemma~5; and Basse and Pedersen \cite{BassePedersen}, Lemma~4.9.
Stochastic Fubini-type results for semimartingales can be founded in
Protter \cite{Protter} and Ikeda and Watanabe \cite{IandW}; however, the
assumptions in these results
are too strong for our purpose.

\begin{theorem}[(Fubini)]\label{stoc_fubini}
Let $f\dvtx T\times S\mapsto\R$ be an
$\mathcal B(T)\otimes\sigma(\mathcal S)$-measurable function such that
%
\begin{equation}
\label{eq:21}
f_x=f(x,\cdot)\in L^\phi\qquad \mbox{for } x\in T \mbox{ and } \int_E
\|f_x\|_\phi \mu(\d{x})<\infty.
\end{equation}
Then $f(\cdot,s)\in L^1(\mu)$ for $m$-a.a. $s\in S$ and $s\mapsto
\int_T
f(x,s)\mu(\d{x})$ belongs to $L^\phi$, all of the below integrals
exist and
%
\begin{equation}
\label{eq:15}
\int_T \biggl(\int_S f(x,s)\Lambda(\d{s})\biggr)\mu(\d{x})=\int_S \biggl(\int_T
f(x,s)\mu(\d{x}) \biggr)\Lambda(\d{s}) \qquad\mbox{a.s.}
\end{equation}
\end{theorem}

\begin{remark}
If $\mu$ is a finite measure, then the last
condition in \eqref{eq:21} is equivalent to
\[
\int_T\biggl[\int_S f(x,s)^2\sigma^2(s)+\int_\R\bigl(|uf(x,s)|^2\wedge|uf(x,s)|\bigr)\nu(\d{u},s)\biggr]m(\d{s})\mu(\d{x})<\infty.
\]
\end{remark}

We will need Theorem~\ref{stoc_fubini} to be able to prove
Proposition~\ref{moving_average_rep}. That proposition yields, in
particular, examples for which the conditions of Theorem 3.1 are
fulfilled. But before proving Theorem~\ref{stoc_fubini}, we will need
the following observation.

\begin{lemma}\label{remark_int}
For all measurable functions $\morf{f}{T\times S}{\R}$ we have
%
\begin{equation}
\label{eq:101}
\biggl\|\int_T |f(x,\cdot)|\mu(\d{x})\biggr\|_\phi\leq\int_T
\norm{f(x,\cdot)}_\phi\mu(\d{x}).
\end{equation}
Moreover, if $\morf{f}{T\times S}{\R}$ is a measurable function such
that $\int_T \norm{f(x,\cdot)}_\phi\mu(\d{x})<\infty$, then for
$m$-a.a. $s\in S$, $f(\cdot,s)\in L^1(\mu)$ and $s\mapsto\int_T
f(x,s)\mu(\d{x})$ is a well-defined function that belongs to
$L^\phi$.
\end{lemma}

\begin{pf}
Let us sketch the proof of \eqref{eq:101}. For $f$ of the form
\[
f(x,s)=\sum_{i=1}^k g_i(s) 1_{A_i}(x),
\]
where $k\geq1$, $g_1,\dots,g_k\in
L^\phi$ and $A_1,\dots, A_k$ are disjoint measurable subsets of $T$ of
finite \mbox{$\mu$-measure}, \eqref{eq:101} easily follows. Hence, by a
monotone class lemma argument, it is possible to show~\eqref{eq:101}
for all measurable $f$. The second statement is a consequence of~\eqref{eq:101}.
\end{pf}

Recall that if $(F,\norm{\cdot})$ is a separable Banach space,
$\mu$ is a measure on $T$ and $\morf{f}{T}{F}$ is a measurable map
such that $\int_T \|f(x)\|\mu(\d{x})<\infty$, then the Bochner
integral $\mathrm{B}\int_T f(x)\mu(\d{x})$ exists in $F$
and $\|\mathrm{B}\int_T f(x)\mu(\d{x})\| \leq \int_T \|f(x)\|\mu
(\d{x})$. Even though
$(L^\phi,\norm{\cdot}_\phi)$ is a Banach space, this result does
not cover Lemma~\ref{remark_int}.

\begin{pf*}{Proof of Theorem~\ref{stoc_fubini}}
For $f$ of the form
%
\begin{equation}
\label{eq:16}
f(x,s)=\sum_{i=1}^n \alpha_i 1_{A_i}(x)1_{B_i}(s),\qquad x\in T,
s\in S,
\end{equation}
where $n\geq1$, $A_1,\dots,A_n$ are measurable
subsets of $T$ of finite $\mu$-measure, $B_1,\dots,B_n\in\mathcal
S$ and $\alpha_1,\dots,\alpha_n\in\R$, the theorem is
trivially true. Thus, for a general $f$ as in the theorem, choose
$f_n$ for $n\geq1$ of the form \eqref{eq:16} such that
$\int_T\|f_n(x,\cdot)-f(x,\cdot)\|_\phi\mu(\d{x})\rightarrow
0$. Indeed, the existence of such a sequence follows by an
application of the monotone class lemma. Let
\[
X_n=\int_E \biggl(\int_S f_n(x,s)\Lambda(\d{s})\biggr)\mu(\d{x}),\qquad
X=\int_E \biggl(\int_S f(x,s)\Lambda(\d{s})\biggr)\mu(\d{x}),
\]
and let us show that $X$ is well defined and
$X_n\rightarrow X$ in $L^1$. This follows since
\[
\mathrm{E}\biggl[\int_E \biggl|\int_S
f(x,s)\Lambda(\d{s})\biggr|\mu(\d{x})\biggr]\leq 2c_2 \int_E
\|f(x,\cdot)\|_\phi\mu(\d{x})<\infty
\]
and
\[
\mathrm{E}[|X_n-X|]\leq 2 c_2 \int_E \|f_n(x,\cdot)-f(x,\cdot)\|_\phi
\mu(\d{x}).
\]
Similarly, let
\[
Y_n=\int_S \biggl(\int_E
f_n(x,s)\mu(\d{x}) \biggr)\Lambda(\d{s}),\qquad Y=\int_S \biggl(\int_E
f(x,s)\mu(\d{x}) \biggr)\Lambda(\d{s})
\]
and let us show that $Y$ is well defined and $Y_n\rightarrow Y$
in $L^1$. By Remark~\ref{remark_int}, $s\mapsto\int_E
f(x,s)\mu(\d{x})$ is a well-defined function that belongs to
$L^\phi$, which shows that $Y$ is well defined.
By \eqref{eq:25} and \eqref{eq:101} we have
\[
\mathrm{E}[|Y_n-Y|]\leq2 c_2 \int_E
\|f_n(x,\cdot)-f(x,\cdot)\|_\phi\mu(\d{x}),
\]
which shows that $Y_n\rightarrow Y$ in $L^1$.
We have, therefore, proved \eqref{eq:15}, since $Y_n=X_n$ a.s.,
$X_n\rightarrow X$ and $Y_n\rightarrow Y$ in $L^1$.
\end{pf*}

Let $Z=(Z_t)_{t\in\R}$ denote an integrable and centered L\'{e}vy
process with L\'{e}vy measure $\nu$ and Gaussian component $\sigma^2$.
Then $Z$ induces a L\'{e}vy basis $\Lambda$ on $S=\R$ and $\mathcal
S=\mathcal B_b(\R)$, the bounded Borel sets, which is uniquely
determined by $\Lambda((a,b])=Z_b-Z_a$ for all $a,b\in\R$ with $a<b$.
In this case $m$ is the Lebesgue measure on $\R$ and
\[
\phi(y,s)=\phi(y)=\sigma^2+\int_\R\bigl(|uy|^2 1_{\{|uy|\leq
1\}}+(2|uy|-1)1_{\{|uy|>1\}}\bigr)\nu(\d{u}).
\]
We will write $\int f(s)\d{Z} _s$ instead of $\int
f(s)\Lambda(\d{s})$. Note that,
$\int_\R f(s)\,\d{Z} _s$ exists and is
integrable if and only if $f\in L^\phi$, that is,
%
\begin{equation}
\label{eq:52}
\int_\R
\biggl(f(s)^2\sigma^2+\int_\R
\bigl(|uf(s)|^2\wedge|uf(s)|\bigr)\nu(\d{x})\biggr)\d{s}<\infty.
\end{equation}
Moreover, if $Z$ is a symmetric $\alpha$-stable L\'{e}vy process,
$\alpha\in(0,2]$, then $L^\phi=L^\alpha(\R,\lambda)$, where
$L^\alpha(\R,\lambda)$ is the space of $\alpha$-integrable functions
with respect to the Lebesgue measure $\lambda$.

\section{Moving average representations}\label{ma}

In wide generality, if $X$ is a continuous-time stationary process, then it is
representable, in principle, as a moving average (MA), that is,
\[
X_{t}=\int_{-\infty}^{t}\psi(t-s)\, \mathrm{d}\Xi_{s},
\]
where $\psi$ is a deterministic function and $\Xi$ has stationary
and orthogonal increments, at least in the second-order sense. (For a
precise statement, see the beginning of Section~\ref{W_K}.) However, an
explicit expression for $\phi$ is seldom available.

We show in Section~\ref{e_MA_rep} that an expression can be found in
cases where the process $X$ is the stationary solution to a Langevin
equation for which the driving noise process $N$ is a PMA, that~is,
%
\begin{equation}
\label{eq:9}
N_t=\int_\R\bigl(f(t-s)-f(-s)\bigr)\,\d{Z} _s,\qquad t\in\R,
\end{equation}
where $Z=(Z_t)_{t\in\R}$ is a suitable process
specified later on and $\morf{f}{\R}{\R}$ is a deterministic function
for which the integrals exist.

In Section~\ref{a_beha}, continuing the discussion from
Section~\ref{asy_auto},
we use the MA representation to study the asymptotic
behavior of the associated autocovariance functions.
Section~\ref{cancel} comments on a notable cancellation effect. But
first, in Section~4.1 we summarize known results concerning
Wold--Karhunen-type representations of stationary continuous-time
processes.

\subsection{Wold--Karhunen-type decompositions}\label{W_K}

Let $X=(X_t)_{t\in\R}$ be a second-order stationary process of mean
zero and
continuous in quadratic mean. Let $F_X$ denote the spectral measure
of $X$, that is, $F_X$ is a finite and symmetric measure on $\R$ satisfying
\[
\mathrm{E}[X_t X_u]=\int_\R \e^{\mathrm{i}(t-u)x}F_X(\d{x}),\qquad t,u\in\R,
\]
and let $F_X'$ denote the density of the
absolutely continuous part of $F_X$.
For each $t\in\R$ let $\mathcal X_t=\spn\{X_s\dvt s\leq t\}$, $\mathcal
X_{-\infty}=\bigcap_{t\in\R} \mathcal X_t$ and $\mathcal
X_\infty=\spn\{X_s\dvt s\in\R\}$ ($\spn$ denotes the $L^2$-closure of the
linear span). Then $X$ is called deterministic if $\mathcal
X_{-\infty}=\mathcal X_\infty$ and purely non-deterministic if
$\mathcal X_{-\infty}=\{0\}$. The following result, which is due to
Satz~5--6 in Karhunen \cite{Karhunen} (cf. also Doob \cite{Doob}, Chapter~XII,
Theorem~5.3), provides a decomposition of stationary processes as a sum
of a deterministic process and a purely non-deterministic process.

\begin{theorem}[(Karhunen)]
Let $X$ and $F_X$ be given as above. If
%
\begin{equation}
\label{eq:28}
\int_\R\frac{\abs{\log F_X'(x)}}{1+x^2}\,\d{x}<\infty,
\end{equation}
then there
exists a unique decomposition of $X$ as
%
\begin{equation}
\label{eq:13}
X_t=\int_{-\infty}^t \psi(t-s)\,\d\Xi_s+V_t,\qquad t\in\R,
\end{equation}
where $\morf{\psi}{\R}{\R}$ is a Lebesgue
square-integrable deterministic function and $\Xi$ is a process with
second-order stationary and orthogonal increments,
$\mathrm{E}[\abs{\Xi_u-\Xi_s}^2]=\abs{u-s}$. Furthermore, for all $t\in\R$,
$\mathcal X_t=\spn\{\Xi_s-\Xi_u\dvt -\infty<u<s\leq t\}$, and $V$
is a deterministic second-order stationary
process.

Moreover, if $F_X$ is absolutely continuous and \eqref{eq:28} is
satisfied, then $V\equiv0$ and hence $X$ is a
backward MA.
Finally, the integral in \eqref{eq:28} is infinite if and only
if $X$ is deterministic.
\end{theorem}

The results in Karhunen \cite{Karhunen} are formulated for complex-valued
processes; however, if~$X$ is real-valued (as it is in our case), then
one can show that all the above processes and functions are real-valued
as well. Note also that if $X$ is Gaussian, then the process $\Xi$ in
\eqref{eq:13} is a standard Brownian motion. If $\sigma$ is a
stationary process with $\mathrm{E}[\sigma^2_0]=1$ and~$B$ is a Brownian
motion, then $\mathrm{d}\Xi_{s}=\sigma_{s}\mathrm{d}B_{s}$ is of the
above type.

A generalization of the classical Wold--Karhunen result to a broad
range of non-Gaussian, infinitely divisible processes was given in
Rosi{\'n}ski \cite{Ros07}.

\subsection{Explicit MA solutions of Langevin equations} \label{e_MA_rep}

Assume initially that $Z$ is an integrable and centered L\'{e}vy
process, and recall that $L^\phi$ is the space of all measurable
functions $\morf{f}{\R}{\R}$ satisfying \eqref{eq:52}. Let
$f\dvtx \R\rightarrow\R$ be a~measurable function
such that $f(t-\cdot)-f(-\cdot)\in
L^\phi$ for all $t\in\R$, and let $N$ be given~by
\[
N_t=\int_\R\bigl(f(t-s)-f(-s)\bigr)\,\d{Z} _s,\qquad t\in\R.
\]

\begin{proposition}\label{moving_average_rep}
Let $N$ be given as above. Then there exists a unique-in-law QOU
process $X$ driven by $N$ with parameter $\lambda>0$, and $X$ is
an MA of the form
%
\begin{equation}
\label{eq:34}
X_t= \int_\R\psi_f(t-s)\,\d{Z} _s, \qquad t\in\R,
\end{equation}
where $\morf{\psi_f}{\R}{\R}$ belongs to $L^\phi$ and
is given by
\begin{equation}
\label{def_of_psi}
\psi_f(t)= \biggl(f(t)-\lambda \e^{-\lambda t} \int_{-\infty}^t
\e^{\lambda
s} f(s)\,\d{s}\biggr), \qquad t\in\R.
\end{equation}
\end{proposition}

\begin{pf}
Since $(t,s)\mapsto f(t-s)-f(-s)$ is measurable we may choose a
measurable modification of $N$ -- see Section~\ref{fubini} --
and hence, by Theorem~\ref{lemma_Masani}, there exists a unique-in-law
QOU process $X$ driven by $N$ with parameter $\lambda$.
For fixed $t\in\R$, we have by \eqref{eq:88} and with
$h_u(s)=f(t-s)-f(t+u-s)$ for
all $u,s\in\R$ and $\mu(\d{u})=1_{\{u\leq
0\}}\e^{\lambda u}\,\d{u}$ that
\[
X_t=\lambda\int_{-\infty}^0 \e^{\lambda u}
(N_t-N_{t+u})\,\d{u}=\int_{-\infty}^0 \biggl(\int_\R h_u(s)\,\d{Z}
_s\biggr)\mu(\d{u}).
\]
By Theorem~\ref{thm_cont} there exist
$\alpha,\beta>0$ such that $\norm{h_u}_\phi\leq\alpha+\beta\abs{t}$
for all $u\in\R$, implying that $\int_\R
\norm{h_u}_\phi\mu(\d{u})<\infty$. By Theorem~\ref{stoc_fubini},
$(u\mapsto h_u(s))\in L^1(\mu)$ for Lebesgue almost all $s\in\R$, which
implies that $\int_{-\infty}^t \abs{f(u)}\e^{\lambda u}\,\d
{u}<\infty$
for all $t>0$, and hence $\psi_f$, defined in~\eqref{def_of_psi}, is a
well-defined function. Moreover, by Theorem~\ref{stoc_fubini},
$\psi_f\in L^\phi(\R,\lambda)$ and
\[
X_t= \int_\R\biggl(\int_{-\infty}^0
h(u,s)\mu(\d{u})\biggr)\,\d{Z} _s=\int_\R\psi_f(t-s)\,\d{Z} _s,\qquad
t\in\R,
\]
which completes the proof.
\end{pf}

Note that for $f=1_{\R_+}$, we have $N_t=Z_t$ and
$\psi_f(t)=\e^{-\lambda t}1_{\R_+}(t)$. Thus, in this case we recover
the well-known result that the QOU process $X$ driven by $Z$ with
parameter
$\lambda>0$ is an MA of the
form $X_t=\int_{-\infty}^t \e^{-\lambda(t-s)} \,\d{Z} _s$.

Let us use the notation $x_+:=x 1_{\{ x\geq0\}}$, and let $c_H$ be
given by
\[
c_H=\frac{\sqrt{2H\sin(\uppi H)\Gamma(2H)}}{\Gamma(H+1/2)}.
\]
A PMA $N$ of the form \eqref{eq:9},
where $Z$ is an $\alpha$-stable L\'{e}vy process with $\alpha\in
(0,2]$ and $f$ is given by $t\mapsto c_H t_+^{H-1/\alpha}$, is called
a \textit{linear fractional $\alpha$-stable motion} of index $H\in
(0,1)$; see Samorodnitsky and Taqqu \cite{Stable}. Moreover, PMAs with $f(t)=t^{\alpha}$ for
$\alpha\in(0,\frac{1}{2})$ and where $Z$ is a square-integrable and
centered L\'{e}vy process are called \textit{fractional L\'{e}vy
processes} in Marquardt \cite{Tina}; these processes provide examples of $f$ and
$Z$ for which Proposition 4.2 applies. Moreover, \cite{Tina},
Theorems~6.2 and 6.3, studies MAs driven by fractional L\'{e}vy processes,
which in some cases also have a representation of the form~%
\eqref{eq:34}.

\begin{corollary}\label{ma_frac}
Let $\alpha\in(1,2]$ and $N$ be a linear fractional
$\alpha$-stable motion of index $H\in(0,1)$. Then there exists
a unique-in-law QOU process $X$ driven by $N$ with parameter
$\lambda>0$, and $X$ is an MA of the form
\[
X_t = \int_{-\infty}^t \psi_{\alpha,H}(t-s)\,\d{Z} _s,
\qquad t\in\R,
\]
where $\psi_{\alpha,H}\dvtx \R_+\rightarrow\R$ is given
by
\[
\psi_{\alpha,H}(t)= c_H\biggl(t^{H-1/\alpha}-\lambda
\e^{-\lambda t}\int_0^t \e^{\lambda
u}u^{H-1/\alpha}\,\d{u}\biggr),\qquad t\geq0.
\]
For $t\to\infty$, we have $\psi_{\alpha,H}(t)\sim
(c_H(H-1/\alpha)/\lambda)t^{H-1/\alpha-1}$, and for $t\to0$,
$\psi_{\alpha,H}(t)\sim c_H t^{H-1/\alpha}$.
\end{corollary}

\begin{remark}
A QOU process driven by a linear fractional $\alpha$-stable motion is
called a~fractional Ornstein--Uhlenbeck process. In
Maejima and Yamamoto \cite{MaejimafsOU}, the existence of the fractional
Ornstein--Uhlenbeck process is shown in the case where $\alpha>1$ and
$1/\alpha<H<1$. (The case $H=1/\alpha$ is trivial since $X=N$.)
The existence in the case $H\in(0,1/\alpha)$ (see Corollary~\ref{ma_frac})
is somewhat unexpected due to
the fact that the sample paths of the linear fractional $\alpha$-stable
motion are unbounded on each compact interval; cf. page~4 in
Maejima and Yamamoto \cite{MaejimafsOU}, where non-existence is surmised.
In the case $\alpha=2$ (i.e., $N$ is a fractional Brownian motion),
Cheridito \textit{et~al.} \cite{CheKawMae03} show the existence of the fractional
Ornstein--Uhlenbeck process.
\end{remark}

In the next lemma we will show a special property of $\psi_f$, given
by \eqref{def_of_psi}; namely that $\int_0^\infty\psi_f(s)\,\d{s}=0$
whenever this integral
is well defined and $f$ tends to zero at $\infty$. This property has a
great impact on the behavior of the autocovariance function of QOU
processes. We
will return to this point in Section~\ref{cancel}.

\begin{lemma}
\label{remark_f}
Let $\morf{f}{\R}{\R}$ be a locally integrable
function that is zero on $(-\infty,0)$ and $\lim_{t\to\infty}f(t)=0$.
Then, $\lim_{t\to\infty}\int_0^t \psi_f(s)\,\d{s}=0$.
\end{lemma}

\begin{pf}
For $t>0$,
\begin{eqnarray*}
\int_0^t \biggl(\lambda \e^{-\lambda s}\int_0^s \e^{\lambda u}f(u)\,\d{u}\biggr) \,\d{s}
&=&
\int_0^t \biggl(\int_u^t\lambda\e^{-\lambda s} \,\d{s}\biggr)\e^{\lambda u}f(u)
\,\d{u}\\
&=&
\int_0^t f(u)\,\d{u}-\e^{-\lambda t}\int_0^t \e^{\lambda u} f(u)\,\d{u},
\end{eqnarray*}
and hence, by using that $\lim_{t\to\infty}f(t)=0,$ we obtain
that
\[
\lim_{t\to\infty} \int_0^t \psi_f(s)\,\d{s} =\lim_{t\to\infty} \biggl(
\e^{-\lambda t}
\int_0^t \e^{\lambda u} f(u)\,\d{u}\biggr)=0.
\]
\upqed\end{pf}

Proposition~\ref{moving_average_rep}
carries over to a much more general setting. For example, if
$N$ is of the form
\[
N_{t}=\int_{\R\times V}[ f( t-s,x) -f(-s,x) ] \Lambda(
\mathrm{d}s,\mathrm{d}x),\qquad t\in\R,
\]
where $\Lambda$ is a centered L\'{e}vy basis on $\R\times V$ ($V$ is
a non-empty space) with control measure
$m(\d{s},\d{x})=\mathrm{d}sn(\mathrm{d}x)$; $a(s,x),\sigma^2(s,x)$
and $\nu(\d{u},(s,x))$, from \eqref{eq:36}, do not depend on
$s\in\R$; and $f(t-\cdot,\cdot)-f(-\cdot,\cdot)\in L^\phi$ for all
$t\in\R$, then, using Theorems~\ref{thm_cont}, \ref{lemma_Masani} and
\ref{stoc_fubini}, the arguments from
Proposition~\ref{moving_average_rep} show that there exists a
unique-in-law QOU process~$X$ driven by $N$ with parameter $\lambda>0$,
and $X$ is given by
\[
X_t=\int_{\R\times V} \psi_f(t-s,x)\Lambda(\d{s},\d{x}), \qquad t\in\R,
\]
where
\[
\psi_f(s,x)= f(s,x)-\lambda \e^{-\lambda s}\int_{-\infty}^s f(u,x)\e^{\lambda u}\,\d{u}, \qquad s\in\R,x\in V.
\]
We recover Proposition~\ref{moving_average_rep} when $V=\{0\}$ and
$n=\delta_0$ is the Dirac delta measure at 0.

\subsection{Asymptotic behavior of the autocovariance function}\label{a_beha}

The representation, from the previous section, of QOU processes as MAs
enables us to handle the autocovariance function analytically. In
Section~\ref{gen_MA} we discuss how the tail behavior of the kernel
$\psi$ of a general MA process determines that of the covariance
function. By use of those results, Section~\ref{cov_MA} relates the
asymptotic behavior of the kernel of the noise $N$ to the asymptotic
behavior of the autocovariance function of the QOU process $X$ driven
by $N$, both for $t\to0$ and $t\to\infty$.

\subsubsection{Autocovariance function of general MAs} \label{gen_MA}

Let $\psi$ be a Lebesgue square-integrable function and $Z$ be a
centered process with stationary and orthogonal increments. Assume for
simplicity that $Z_0=0$ a.s. and $\var{Z}(t)=t$. Let
$X=\psi*Z=(\int_{-\infty}^t \psi(t-s)\,\d{Z} _s)_{t\in\R}$ be a
backward MA; $\ac{X}$ be its autocovariance function, that is
\[
\ac{X}(t)=\mathrm{E}[X_t X_0]=\int_0^\infty\psi(t+s)\psi(s)\,\d{s},\qquad t\in\R;
\]
and $\acc{X}(t)=\ac{X}(0)-\ac{X}(t)=\frac{1}{2}\mathrm{E}[(X_t-X_0)^2]$.
The behavior of $\ac{X}$ at $0$ or $\infty$ corresponds in large
extent to the behavior of the kernel $\psi$ at $0$ or $\infty$,
respectively.

Indeed, we have the following result, in which
$k_\alpha$ and $j_\alpha$ are constants given by
\begin{eqnarray*}
k_{\alpha}
&=&
\Gamma(1+\alpha)\Gamma(-1-2\alpha)\Gamma(-\alpha)^{-1} ,\qquad \alpha\in(-1,-1/2),
\\
j_\alpha
&=&
(2\alpha+1)\sin\bigl(\uppi(\alpha+1/2)\bigr)\Gamma(2\alpha+1)\Gamma(\alpha+1)^{-2}, \qquad \alpha\in
(-1/2,1/2) .
\end{eqnarray*}

\begin{proposition}\label{pro_cov_MA}
Let the setting be as described
above.
\begin{longlist}[(iii)]
\item\hypertarget{lemma_1}For $t\rightarrow\infty$ and $\alpha\in
(-1,-\frac{1}{2})$, $\psi(t)\sim c t^\alpha$ implies
$\ac{X}(t)\sim(c^2 k_{\alpha})t^{2\alpha+1}$, provided
$|\psi(t)|\leq c_1 t^{\alpha}$ for all $t>0$ and some $c_1>0$.
\item\hypertarget{lemma_2} For $t\rightarrow
\infty$ and $\alpha\in(-\infty,-1)$,
$\psi(t)\sim c t^{\alpha}$ implies $\ac{X}(t)/t^\alpha\to
c \int_0^\infty\psi(s)\,\d{s}$, and hence $\ac{X}(t)\sim
(c \int_0^\infty\psi(s)\,\d{s})t^\alpha$, provided $\int_0^\infty
\psi(s)\,\d{s}\neq0$.
\item\hypertarget{lemma_3}For $t\rightarrow0$ and $\alpha\in
(-\frac{1}{2},\frac{1}{2})$,
$\psi(t)\sim c t^\alpha$
implies $\acc{X}(t)\sim(c^2 j_{\alpha}/2)\abs{t}^{2\alpha+1}$,
provided $\psi$ is absolutely continuous on $(0,\infty)$ with density
$\psi'$ satisfying $|\psi'(t)|\leq c_2 t^{\alpha-1}$ for all $t>0$ and
some $c_2>0$.
\end{longlist}
\end{proposition}

\begin{pf}
\hyperlink{lemma_1}{(i)} Let $
\alpha\in(-1,-\frac{1}{2})$ and assume that $\psi(t)\sim c t^\alpha$
as $t\rightarrow\infty$ and $|\psi(t)|\leq c_1 t^{\alpha}$ for $t>0$.
Then
%
\begin{eqnarray}\label{eq:130}
\ac{X}(t)
&=&
\int_0^\infty\psi(t+s)\psi(s)\,\d{s}\nonumber
\\
&=& t\int_0^\infty\psi\bigl(t(s+1)\bigr)\psi(ts)\,\d{s}\nonumber
\\[-8pt]\\[-8pt]
&=&
t^{2\alpha+1} \int_0^\infty\frac{\psi(t(1+s)) \psi(ts)} {(t(1+s))^{\alpha}(ts)^{\alpha}}(1+s)^{\alpha} s^{\alpha} \,\d{s}\nonumber
\\
&\sim&
t^{2\alpha+1} c^2 \int_0^\infty(1+s)^{\alpha} s^{\alpha} \,\d{s} \qquad\mbox{as }t\rightarrow\infty.\nonumber
\end{eqnarray}
Since
\[
\int_0^\infty(1+s)^{\alpha} s^{\alpha}\,\d{s}=\frac{\Gamma(1+\alpha)\Gamma(-1-2\alpha)}{\Gamma(-\alpha)}=k_\alpha,
\]
\eqref{eq:130} shows that
$\ac{X}(t)\sim(c^2 k_{\alpha})t^{2\alpha+1} $ for
$t\rightarrow\infty$.

\hyperlink{lemma_2}{(ii)} Let $\alpha\in(-\infty,-1)$ and assume that
$\psi(t)\sim c t^{\alpha}$ for $t\rightarrow\infty$. Note that
$\psi\in
L^1(\R_+,\lambda)$ and for some $K>0$ we have for all $t\geq K$ and
$s>0$ that $|\psi(t+s)|/t^\alpha\leq2|c|(t+s)^\alpha/t^\alpha\leq
2|c|$. Hence, by applying Lebesgue's dominated convergence theorem, we
obtain
\[
\ac{X}(t)=t^\alpha\int_0^\infty
\biggl(\frac{\psi(t+s)}{t^\alpha}\psi(s)\biggr)\,\d{s} \sim t^\alpha c \int
_0^\infty
\psi(s)\,\d{s} \qquad\mbox{for }t\rightarrow\infty.
\]

\hyperlink{lemma_3}{(iii)} By letting
\[
f_t(s):= \frac{\psi(t(s+1))-\psi(ts)}{t^\alpha},\qquad t>0, s\in\R,
\]
we have
\begin{equation}
\label{eq:1}
\mathrm{E}[(X_t-X_0)^2]=  t\int\bigl[\psi\bigl(t(s+1)\bigr)-\psi(ts)\bigr]^2\,\d{s}=t^{2\alpha+1} \int \vert f_t(s)\vert^2\,\d{s}.
\end{equation}
As $t\rightarrow0$, we find
\[
f_t(s)=\frac{\psi(t(s+1))}{(t(s+1))^\alpha}
(s+1)^\alpha
-\frac{\psi(ts)}{(ts)^\alpha}s^\alpha\rightarrow c\bigl((s+1)^\alpha_+
-s^\alpha_+\bigr).
\]
Choose $\delta>0$ such that $|\psi(x)|\leq2 x^\alpha$ for $x\in
(0,\delta)$. By our assumptions we have for all $s\geq\delta$ that
\begin{eqnarray*}
\abs{f_t(s)}
&=&
t^{-\alpha} \biggl\vert\int_{ts}^{t(1+s)} \psi'(u)\,\d{u}\biggr\vert\leq t^{-\alpha+1}\sup_{u\in[st,t(s+1)]}\abs{\psi'(u)}
\\
&\leq&
c_2 t^{-\alpha+1}\sup_{u\in[st,t(s+1)]} \abs{u}^{\alpha-1}= c_2 t^{-\alpha+1}\abs{ts}^{\alpha-1} = c_2s^{\alpha-1},
\end{eqnarray*}
and for $s\in[-1,\delta)$,
$|f_t(s)|\leq 2c[ (1+s)^\alpha+s^\alpha_+]$.
This shows that there exists a function $g\in L^2(\R_+,\lambda)$ such
that
$\abs{f_t}\leq g$ for all $t>0$, and thus, by Lebesgue's dominated
convergence theorem, we have
%
\begin{equation}
\label{eq:41}
\int \vert f_t(s)\vert^2\,\d{s}\mathop{\longrightarrow}_{t\to0}
c^2\int \bigl((s+1)^\alpha_+
-s^\alpha_+ \bigr)^2\,\d{s}=c^2 j_{\alpha}.
\end{equation}
Together with \eqref{eq:1}, \eqref{eq:41} shows that
$\acc{X}(t)\sim( c^2
j_\alpha/2)t^{2\alpha+1} $ for $t\rightarrow0$.
\end{pf}

\begin{remark}
It would be of interest to obtain a general result covering
Proposition~\ref{pro_cov_MA}\hyperlink{lemma_2}{(ii)} in the case
$\int_0^\infty\psi(s)\,\d{s}=0$. Recall that $\psi_f$, given by
\eqref{def_of_psi}, often satisfies that $\int_0^\infty
\psi_f(s)\,\d{s}=0$, according to Lemma~\ref{remark_f}.
\end{remark}

\begin{example}
Consider the case where $\psi(t)=t^{\alpha} \e^{-\lambda t}$ for
$\alpha\in(-\frac{1}{2},\infty)$ and $\lambda>0$. For
$t\rightarrow
0$, $\psi(t)\sim t^\alpha$, and hence
$\acc{X}(t)\sim(j_\alpha/2)t^{2\alpha+1}$ for $t\rightarrow0$ and
$\alpha\in(-\frac{1}{2},\frac{1}{2})$, by
Proposition~\ref{pro_cov_MA}\hyperlink{lemma_3}{(iii)}
(compare with Barndorff-Nielsen \textit{et~al}. \cite{BNCorPod09}).
\end{example}

Note that if $X=\psi*Z$ is a moving average, as above, then by
Proposition~\ref{pro_cov_MA}\hyperlink{lemma_1}{(i)} and for $t\to\infty$,
$\ac{X}(t)\sim c_1 t^{-\alpha}$ with $\alpha\in(0,1)$, provided that
$\psi(t)\sim c_2 t^{-(\alpha+1)/2}$ and
$\abs{\psi(t)}\leq c_3 t^{-(\alpha+1)/2}$. This shows that $X$ has
long-range dependence of order $\alpha$.

Let us conclude this section with a short discussion of when an MA
$X=\psi*Z$ is a~semimartingale. It is often very important that the
process of interest is a semimartingale, especially in finance, where
the semimartingale property of the asset price is equivalent to the
property that the capital process depends continuously on the chosen
strategy; see Section~8.1.1 in Cont and Tankov \cite{ContTankov}. In the case where
$Z$ is a Brownian motion, Theorem~6.5, in Knight \cite{Knight} shows that $X$
is an $\F^Z$ semimartingale if and only if $\psi$ is absolutely
continuous on $[0,\infty)$ with a square-integrable density. (Here
$\F^Z_t:=\sigma(Z_s\dvt s\in(-\infty,t])$.) For a further study of the
semimartingale property of PMA and more general processes, see
\cite{Andreas3,Andreas1,Andreas2} in the Gaussian case, and
Basse and Pedersen \cite{BassePedersen} for the infinitely divisible case.

\subsubsection{QOU processes with PMA noise}\label{cov_MA}

Let us return to the case of a QOU process driven by a PMA. Let $Z$ be
a centered L\'{e}vy process, $\morf{f}{\R}{\R}$ be a measurable
function that is 0 on $(-\infty,0)$ and satisfies
$f(t-\cdot)-f(-\cdot)\in L^\phi$ for all $t\in\R$ and $N$ be given by
%
\begin{equation}\label{pma_noise}
N_t=\int_\R[ f(t-s)-f(-s)]\,\d{Z} _s,\qquad t \in\R.
\end{equation}
First, we will consider the relationship between the behavior of the
kernel of the noise $N$ and that of the kernel $\psi_f$ of the
corresponding moving average $X$.

\begin{proposition}\label{cor_MA_noise}
Let $N$ be given by \eqref{pma_noise}, and $X$ be a QOU process driven
by $N$ with parameter $\lambda>0$.
\begin{longlist}[(ii)]
\item\hypertarget{pma_1} Let $\alpha\in(-1,-\frac{1}{2})$ and assume that,
for some $c\neq0$, $f$ is continuous differentiable in a
neighborhood of $\infty$ with $f'(t)\sim c t^\alpha$ for $t\to\infty$.
Then, for $t\to\infty$, we have $\ac{X}(t)\sim
(\frac{c^2k_\alpha}{\lambda^2})t^{2\alpha+1}$, provided $|f(t)|\leq r
t^\alpha$ for all $t>0$ and some $r>0$.
\item\hypertarget{pma_2} Let
$\alpha\in(-\frac{1}{2},\frac{1}{2})$ and $f(t)\sim ct^\alpha$ for
$t\to0$. Then, for $t\to0$, we have $\acc{X}(t)\sim(c^2
j_\alpha/2)\abs{t}^{2\alpha+1}$, provided $f$ is two times continuous
differentiable in a neighborhood of $\infty$ with
$f''(t)=\mathrm{O}(t^{\alpha-1})$ for $t\to\infty$, and that $f$ is absolutely
continuous on $(0,\infty)$ with a density $f'$ satisfying $\sup_{t\in
(0,t_o)}\abs{f'(t)}t^{1-\alpha}<\infty$ for all $t_0>0$.
\end{longlist}
\end{proposition}

\begin{pf}
\hyperlink{pma_1}{(i)} Choose $\beta>0$ such that $f$ is continuous
differentiable on $[\beta,\infty)$. By partial integration, we have for
$t\geq\beta$,
%
\begin{equation}
\label{eq:107}
\psi_f(t)=\e^{-\lambda t}\biggl(\e^{\lambda a }f(a)-\lambda\int_{-\infty}^a
\e^{\lambda
s}f(s)\,\d{s}\biggr)+ \e^{-\lambda t} \int_a^t \e^{\lambda s}
f'(s)\,\d{s},
\end{equation}
showing that $\psi_f(t)\sim(\frac{c}{\lambda})t^\alpha$ for
$t\to\infty$.
Choose $k>0$ such that $\abs{\psi_f(t)}\leq(2c/\lambda)t^\alpha$
for all
$t\geq k$. By \eqref{def_of_psi} we have that $\sup_{t\in
[0,k]}\abs{\psi_f(t)t^{-\alpha}}<\infty,$
since $\sup_{t\in[0,k]}\abs{f(t)t^{-\alpha}}<\infty$, and hence
there exists a constant $c_1>0$ such that
$\abs{\psi_f(t)}\leq c_1 t^\alpha$ for all $t>0$. Therefore,
(i) follows by
Proposition~\ref{pro_cov_MA}\hyperlink{lemma_1}{(i)}.

\hyperlink{pma_2}{(ii)} Choose $\beta>0$ such that $f$ is two times continuous
differentiable on $[\beta,\infty)$. By \eqref{eq:107} and partial
integration we have for $t>\beta$ and $t\to\infty$,
\begin{eqnarray*}
\psi_f'(t)
&=&
f'(t)-\lambda\psi_f(t)=f'(t)-\lambda \e^{-\lambda t}\int_\beta^t \e^{\lambda s}f'(s)\,\d{s}+\mathrm{O}(\e^{-\lambda t})
\\
&=&
\e^{-\lambda t}\int_\beta^t\e^{\lambda s} f''(s)\,\d{s}+\mathrm{O}(\e^{-\lambda t})=\mathrm{O}(t^{\alpha-1}),
\end{eqnarray*}
where we in the last equality have used that
$f''(t)=\mathrm{O}(t^{\alpha-1})$ for $t\to\infty$. Using that
$\abs{\psi_f'(t)}\leq\abs{f'(t)} +\lambda\abs{\psi_f(t)}$ and
$\sup_{t\in
(0,t_0)}\abs{f'(t)t^{1-\alpha}}<\infty$ for all $t_0>0$, it follows
that there exists a $c_2>0$ such that $\abs{\psi_f'(t)}\leq c_1
t^{\alpha-1}$ for all
$t>0$. Moreover,\vspace*{2pt} for $t\to0$, we have that
$\psi_f(t)\sim ct^\alpha$. Hence, (ii) follows by
Proposition~\ref{pro_cov_MA}\hyperlink{lemma_3}{(iii)}.
\end{pf}

Now consider the following set-up: Let $Z=(Z_t)_{t\in\R}$ be a
centered and square-integrable L\'{e}vy process, and for $H\in(0,1)$,
$r_0\neq0$, $\delta\geq0$, let
%
\begin{equation}
\label{eq:108}
f(t)=r_0 (\delta\vee t)^{H-1/2} \quad\mbox{and}\quad
N_t^{H,\delta}=\int_\R[f(t-s)-f(-s)]\,\d{Z}_s.
\end{equation}
Note that when $\delta=0$ and $Z$ is a Brownian motion, $N^{H,\delta}$
is a constant times the fBm of index $H$, and when $\delta>0$,
$N^{H,\delta}$ is a semimartingale. We have the following corollary to
Proposition~\ref{cor_MA_noise}.

\begin{corollary}\label{exp_fBm}
Let $N^{H,\delta}$ be given by \eqref{eq:108}, and let $X^{H,\delta}$
be a QOU process driven by
$N^{H,\delta}$ with parameter $\lambda>0$.
Then, for $H\in(\frac{1}{2},1)$ and $t\to\infty$,
\[
\ac{X^{H,\delta}}(t)\sim \bigl(r_0^2 k_{H-3/2}(H-1/2)/\lambda^2\bigr)t^{2H-2},\qquad \delta\geq0,
\]
and for $H\in(0,1)$ and $t\to0$,
%
\begin{equation}
\label{eq:43}
\acc{X^{H,\delta}}(t)\sim
\cases{
(r_0^2\delta^{2H-1}/2)\abs{t}, &\quad$ \delta>0$,\cr
(r_0^2j_{H-1/2}/2)\abs{t}^{2H}, &\quad$ \delta=0$.
}
\end{equation}
\end{corollary}

\begin{pf}
For $H\in(\frac{1}{2},1)$, let $\beta=\delta$. Then $f\in
C^1((\beta,\infty);\R)$ and, for $t>\beta$, $f'(t)=ct^{\alpha}$, where
$\alpha=H-3/2\in(-1,-\frac{1}{2})$ and $c=r(H-1/2)$. Moreover,
$\abs{f(t)}\leq r\delta t^{\alpha}$. Thus,
Proposition~\ref{cor_MA_noise}\hyperlink{lemma_1}{(i)} shows that
$\ac{X^{H,\delta}}(t)\sim
(c^2k_\alpha/\lambda^2)t^{2\alpha+1}=(r^2(H-1/2)^2k_{H-3/2}/\lambda
^2)t^{2H-2}$.
To show \eqref{eq:43} assume that $H\in(0,1)$. For $t\to0$, we have
$f(t)\sim c t^{\alpha}$, where $c=r_0$ and $\alpha=H-1/2\in
(-\frac{1}{2},\frac{1}{2})$ when $\delta=0$, and $c=r_0 \delta^{H-1/2}$
and $\alpha=0$ when $\delta>0$. For $\beta=\delta$, $f\in
C^2((\beta,\infty);\R)$ with $f''(t)=r_0(H-1/2)(H-3/2)t^{H-5/2}$,
showing that $f''(t)=\mathrm{O}(t^{\alpha-1})$ for $t\to\infty$ (both for
$\delta>0$ and $\delta=0$). Moreover, $f$ is absolutely continuous on
$(0,\infty)$ with density
$f'(t)=r_0(H-1/2)t^{H-3/2}1_{[\delta,\infty)}(t)$. This shows that
$\sup_{t\in(0,t_0)}\abs{f'(t)t^{1-\alpha}}<\infty$ for all
$t_0>0$ (both for $\delta>0$ and $\delta=0$). Hence \eqref{eq:43}
follows by Proposition~\ref{cor_MA_noise}\hyperlink{pma_2}{(ii)}.
\end{pf}

\subsection{Stability of the autocovariance function}
\label{cancel}

Let $N$ be a PMA of the form \eqref{eq:9},
where $Z$ is a centered square-integrable
L\'{e}vy process, and $f(t)=c_H t_+^{H-1/2}$, where $H\in(0,1)$.
(Recall that if $Z$ is a Brownian motion, then $N$ is an fBm of index
$H$.) Let $X$ be a QOU process driven by $N$ with parameter
$\lambda>0$, and recall that by Proposition~\ref{moving_average_rep},
$X$ is an MA of the form
\[
X_t= \int_{-\infty}^t \psi_H(t-s)\,\d{Z} _s,
\qquad t\in\R,
\]
where
\[
\psi_H(t)=c_H\biggl(t^{H-2/2}-\lambda \e^{-\lambda t}\int_0^t\e^{\lambda u}u^{H-1/2}\,\d{u}\biggr), \qquad t\geq0.
\]
Below we will discus some stability properties for the autocovariance
function under minor modification of the kernel function.

For all bounded measurable functions $\morf{f}{\R_+}{\R}$ with compact
support, let
$X^f_t=\int_{-\infty}^t(\psi_H(t-s)-f(t-s))\,\d{Z} _s$. We
will think of $X^f$ as an MA where we have made a~minor change of
$X$'s kernel. Note that if we let $Y^f_t=X_t-X^f_t=\int_{-\infty}^t
f(t-s)\,\d{Z} _s$, then the autocovariance function $\ac{Y^f}(t)$, of
$Y^f$, is zero whenever $t$ is large enough due to the fact that $f$
has compact support.

\begin{corollary} \label{cor_stab}
We have the following two situations in which $c_1,c_2,c_3\neq0$
are non-zero constants.
\begin{longlist}[(ii)]
\item\hypertarget{stab_1} For $H\in(0,\frac{1}{2})$ and $\int_0^\infty
f(s)\,\d{s}\neq0$, we have for $t\to\infty$,
%
\[
\ac{X^f}(t)\sim c_2 \ac{X}(t)t^{1/2-H}\sim c_1 t^{H-3/2}.
\]
\item\hypertarget{stab_2}
For $H\in(\frac{1}{2},1)$, we have for $t\to\infty$,
\[
\ac{X^f}(t)\sim\ac{X}(t)\sim c_3t^{2H-2}.
\]
\end{longlist}
\end{corollary}

Thus for $H\in(0,\frac{1}{2})$, the above shows that the behavior of
the autocovariance function at infinity is changed dramatically by
making a minor change of the kernel. In particular, if $f$ is a
positive function, not the zero function, then $\ac{X^f}(t)$ behaves
as $t^{1/2-H}\ac{X}(t)$ at infinity. On the other hand, when $H\in
(\frac{1}{2}, 1),$ the behavior of the autocovariance function at
infinity does not change if we make a minor change to the kernel. That
is, in this case the autocovariance function has a stability property,
contrary to the case where $H\in(0,\frac{1}{2})$.

\begin{remark}
Note that the dramatic effect appearing from
Corollary~\ref{cor_stab}\hyperlink{stab_1}{(i)} is associated with the fact
that $\int_0^\infty\psi_H(s)\,\d{s}=0$, as shown in Lemma~\ref{remark_f}.
\end{remark}

\begin{pf*}{Proof of Corollary~\ref{cor_stab}}
By Corollary~\ref{ma_frac} we have for $t\to\infty$ that $\psi
_H(t)\sim
ct^\alpha$, where $c=c_H(H-1/2)/\lambda$ and
$\alpha=H-3/2$. To show \hyperlink{stab_1}{(i)}, assume that $H\in
(0,\frac{1}{2})$ and hence
$\alpha\in(-\infty,-1)$. According to Lemma~\ref{remark_f}, we
have that $\int_0^\infty\psi_H(s)\,\d{s}=0$ and hence $\int
_0^\infty
[\psi_H(s)-f(s)]\,\d{s}\neq0,$ since
$\int_0^\infty f(s)\,\d{s}\neq0$ by assumption.
From Proposition~\ref{pro_cov_MA}\hyperlink{lemma_2}{(ii)} and for $t\to
\infty,$ we have that
$\ac{X^f}(t)(t)\sim c_1 t^{2\alpha+1}=c_1 t^{H-3/2}$, where
$c_1=c\int_0^\infty [\psi_H(s)-f(s)]\,\d{s}$. On the
other hand, by Corollary~\ref{cor_var} we have that $\ac{X}(t)\sim
(H(H-1/2)/\lambda^2)t^{2H-2}$ for $t\to\infty$, and hence we have
shown \hyperlink{stab_1}{(i)} with $c_2 =c_1\lambda^2/(H(H-1/2)$.
For $H\in(\frac{1}{2},1)$ we have that $\alpha\in(-1,
-\frac{1}{2})$, and hence \hyperlink{stab_2}{(ii)} follows by
Proposition~\ref{pro_cov_MA}\hyperlink{lemma_1}{(i)}.
\end{pf*}


\begin{appendix}\label{app}
\setcounter{lemma}{0}
\setcounter{equation}{0}
\section*{Appendix}


In this Appendix we will show an auxiliary continuity result used
several times in the paper. The main result in this Appendix is
Theorem~\ref{thm_cont}; Corollary~\ref{p_cont} is used in
Theorem~\ref{lemma_Masani}, while the general modular setting is needed
to prove Proposition~\ref{moving_average_rep}. For the basic
definitions and properties of linear metric spaces, modulars and
$F$-norms, we refer to Rolewicz \cite{Rolewicz}.

Let $(E,\mathcal E,\mu)$ be a $\sigma$-finite measure space, and
$\morf{\phi}{\R}{\R_+}$ an even and continuous function that is
non-decreasing on $\R_+$, with $\phi(0)=0$. Assume there exists a
constant $C>0$ such that $\phi(2x)\leq C\phi(x)$ for all $x\in\R$
(that is, $\phi$ satisfies the $\Delta_2$ condition). Let
$L^0=L^0(E,\mathcal E,\mu)$ denote the space of all measurable
functions from $E$ into $\R$; $\Phi$ denote the modular on $L^0$
given by
\[
\Phi(g)=\int_E \phi(g)\,\d\mu,\qquad g\in L^0;
\]
and $L^\phi=\{g\in L^0:\Phi(g)<\infty\}$ denote the corresponding
modular space. Furthermore, for $g\in L^0$, define
\[
\rho(g)=\inf\{c>0\dvt \Phi(g/c)\leq c\}\quad \mbox{and}\quad
\norm{g}_\phi=\inf\{c>0\dvt \Phi(g/c)\leq1\}.
\]
Then $\rho$ is an $F$-norm on $L^\phi$ and, in particular,
$d_\phi(f,g)=\rho(f-g)$ is an invariant metric on $L^\phi$. Moreover,
when $\phi$ is convex, the \textit{Luxemburg norm}
$\norm{\cdot}_\phi$ is a norm on $L^\phi$; see Khamsi \cite{Modular}.

\begin{theorem}\label{thm_cont}
Let $\morf{f}{\R\times E}{\R}$ denote a
measurable function satisfying $f_t=f(t,\cdot)\in L^\phi$ for all
$t\in\R$, and
%
\begin{equation}
\label{eq:8}
d_\phi(f_{t+u},f_{v+u})=d_\phi(f_{t}, f_{v}) \qquad\mbox{for all
}t,u,v\in\R.
\end{equation}
Then, $(t\in\R)\mapsto(f_t\in L^\phi)$ is continuous. Moreover, if
$\phi$ is convex, then
there exist $\alpha,\beta>0$ such that
$\norm{f_t}_\phi\leq\alpha+\beta\abs{t}$ for all $t\in\R$.
\end{theorem}

To prove Theorem~\ref{thm_cont}, we shall need the following lemma.

\begin{lemma}\label{Borel_m}
Let $\morf{f}{\R\times E}{\R}$ denote a
measurable function, such that $f_t\in L^\phi$ for all
$t\in\R$. Then, $(t\in\R)\mapsto(f_t\in L^\phi)$
is Borel measurable and has a separable range.
\end{lemma}

Recall that $\morf{f}{E}{F}$ has a separable range, if $f(E)$ is a
separable subset of $F$.

\begin{pf*}{Proof of Lemma \ref{Borel_m}}
We will use a monotone class lemma argument to prove this result,
so let $\mathcal{M}_2$ be the set of all functions $f$ for which
Lemma~\ref{Borel_m} holds
and $\mathcal{M}_1$ be the set of all functions $f$ of the form
\[
f_t(s)=\sum_{i=1}^n \alpha_i 1_{A_i}(t) 1_{B_i}(s),
\qquad t\in\R, s\in E,
\]
where, for $n\geq1$, $A_1,\dots,A_n$ are
measurable subsets of $\R$, $B_1,\ldots,B_n$ are measurable subsets of
$E$ of finite $\mu$ measure and $\alpha_1,\dots,\alpha_n\in\R$.
Let us
show that $\Psi_f\dvtx (t\in\R)\mapsto(f_t\in L^\phi)$ is measurable.
Since, for all $g\in L^\phi$, $t\mapsto d_\phi(f_t,g)$ is measurable,
we get that for all $g\in L^\phi$ and $r>0$, $\Psi_f^{-1}(B(g,r))$ is
measurable (we use the notation, $B(g,r)=\{h\in
L^\phi\dvt d_\phi(g,h)<r\}$). Therefore, since $\Psi_f$ has separable
range, it follows that $\Psi_f$ is measurable (recall that the Borel
$\sigma$-field in a separable metric space is generated by the open
balls).
This shows that $\mathcal{M}_1 \subseteq\mathcal{M}_2$. Note that
the set
$\mathrm{b}\mathcal{M}_2$ of bounded elements from $\mathcal M_2$ is a
vector space with $1\in\mathrm{b}\mathcal{M}_2$, and that
$(f_n)_{n\geq1}\subseteq\mathrm{b}\mathcal{M}_2$ with $0\leq
f_n\uparrow f\leq K$ implies that $f\in\mathrm{b}\mathcal{M}_2$.
Moreover, since $\mathcal{M}_1$ is stable under pointwise
multiplication, the monotone class lemma (see \cite{Rogers},
Chapter~II, Theorem~3.2) shows that
\[
\mathrm{bM}\bigl(\mathcal{B}(\R)\times\mathcal{F}\bigr)=
\mathrm{bM}(\sigma(\mathcal{M}_1))\subseteq
\mathrm{b}\mathcal{M}_2.
\]
(For a family of functions $\mathcal M$, $\sigma(\mathcal M)$ denotes
the least $\sigma$-algebra for which all the functions are measurable,
and for each $\sigma$-algebra $\mathcal E$, $\mathrm{bM}(\mathcal E)$
denotes the space of all bounded $\mathcal E$-measurable functions.)
For a general function $f$, define $f^{(n)}$ by
$f^{(n)}_t=f_t1_{\{|f_t|\leq n\}}$. For all $n\geq1$,
$f^{(n)}$ is a bounded measurable function and hence
$\Psi_{f^{(n)}}$ is a measurable map with a separable range.
Moreover, $\lim_n \Psi_{f^{(n)}}=\Psi_{f}$ pointwise in $L^\phi$,
showing that $\Psi_f$ is measurable and has a separable range.
\end{pf*}

\begin{pf*}{Proof of Theorem~\ref{thm_cont}}
Let $\Psi_f$ denote the map $(t\in\R)\mapsto(f_t\in L^\phi)$,
and for fixed $\epsilon>0$ and arbitrary $t\in\R$, consider the
ball $B_t=\{s\in\R\dvt d_\phi(f_t,f_s)<\epsilon\}$.
By Lemma~\ref{Borel_m}, $\Psi_f$ is
measurable, and hence $B_t$ is a measurable subset of $\R$ for all
$t\in\R$. According to Lemma~\ref{Borel_m}, $\Psi_f$ has a separable
range and, therefore, there exists a countable set $(t_n)_{n\geq
1}\subseteq\R$ such that the range of $\Psi_f$ is included in
$\bigcup_{n\geq
1}B(f_{t_n},\epsilon)$, implying that
$\R=\bigcup_{n\geq 1} B_{t_n}$. In particular, there exists an $n\geq
1$ such that $B_{t_n}$ has a strictly positive Lebesgue measure. By the
Steinhaus lemma, see \cite{Bingham}, Theorem~1.1.1, there exists a
$\delta>0$ such that $(-\delta,\delta)\subseteq B_{t_n}-B_{t_n}$. Note
that by \eqref{eq:8} it is enough to show continuity of $\Psi_f$ at
$t=0$. For $|t|<\delta$ there exists, by definition, $s_1,s_2\in\R$
such that $d_\phi(f_{t_n},f_{s_i})<\epsilon$ for $i=1,2$, showing that
\[
d_\phi(f_t,f_0)\leq d_\phi(f_t,f_{s_1})+d_\phi
(f_t,f_{s_2})<2\epsilon,
\]
which completes the proof of the continuity part.

To show the last part of the theorem, assume that $\phi$ is convex.
For each $t>0$ choose $n=0,1,2,\dots$ such that $n\leq t< n+1$. Then,
\begin{eqnarray}\label{eq:37}
\norm{f_t-f_0}_\phi&
\leq&
\sum_{i=1}^n \norm{f_i-f_{i-1}}_\phi+\norm{f_t-f_n}_\phi\nonumber
\\[-8pt]\\[-8pt]
&\leq&
n \norm{f_1-f_{0}}_\phi+\norm{f_{t-n}-f_0}_\phi\leq t \beta+a,\nonumber
\end{eqnarray}
where $\beta=\norm{f_1-f_0}_\phi$ and $a=\sup_{s\in
[0,1]}\norm{f_s-f_0}_\phi$. We have already shown that
$t\mapsto f_t$ is continuous, and hence $a<\infty$. Since
$\norm{f_{-t}-f_0}_\phi=\norm{f_t-f_0}_\phi$ for all $t\in\R$,
\eqref{eq:37} shows that $\norm{f_t-f_0}_\phi\leq a+\beta\abs{t}$ for
all $t\in\R$, implying that $\norm{f_t}_\phi\leq\alpha+\beta
\abs{t},$ where $\alpha=a+\norm{f_0}_\phi$.
\end{pf*}

For $(E,\mathcal E,\mu)=(\Omega,\F,P)$ and $\phi(t)=\abs{t}^p$ for
$p>0$ or $\phi(t)=\abs{t}\wedge1$ for $p=0$, we have the following
corollary to Theorem~\ref{thm_cont}.

\begin{corollary}\label{p_cont}
Let $p\geq0$ and $X=(X_t)_{t\in\R}$ be a measurable process
with stationary increments and finite $p$ moments. Then $X$ is
continuous in $L^p$. Moreover, if $p\geq1$, then there exist
$\alpha,\beta>0$ such that
$\norm{X_t}_p\leq\alpha+\beta\abs{t}$ for all $t\in\R$.
\end{corollary}

Note that in Corollary~\ref{p_cont} the reversed implication is also
true; in fact, all stochastic processes $X=(X_t)_{t\in\R}$ that are
continuous in $L^0$ have a measurable modification according to
Theorem~2 in Cohn \cite{Cohn}. The idea of using the Steinhaus lemma to
prove Theorem~\ref{thm_cont} is borrowed from
Surgailis \textit{et~al}. \cite{Rosinskipreprint}, where Corollary~\ref{p_cont} is shown for
$p=0$. Furthermore, when $\mu$ is a probability measure and
$\phi(t)=\abs{t}\wedge1$, Lemma~\ref{Borel_m} is known from
Cohn \cite{Cohn}.
\end{appendix}

%

%
\printhistory

\end{document}